\documentclass[a4paper,11pt]{article}

\usepackage{marginnote}
\usepackage{color}

\usepackage{palatino}
\usepackage{mathpazo}

\usepackage{mathrsfs}
\usepackage{amssymb,amsmath,amsthm,url}
\usepackage{verbatim}
\usepackage[small]{titlesec}

\usepackage{enumitem}

\theoremstyle{plain} 
\newtheorem*{thm}{Main Theorem}

\newtheorem{lemma}{Lemma}[section]

\theoremstyle{definition}

\theoremstyle{remark}
\newtheorem{remark}[lemma]{Remark}
\newtheorem*{runex}{Running Example}

\newcommand{\IG}{\textsf{IG}}
\newcommand{\RIG}{\textsf{RIG}}

\newcommand{\st}{\::\:}
\newcommand{\rrestriction}{\!\!\restriction}

\newcommand{\lexless}{<_{\textup{lex}}}
\newcommand{\lexgreater}{>_{\textup{lex}}}
\newcommand{\mmin}{\textup{min}}

\newcommand{\im}{\textup{im}\:\:}
\newcommand{\rank}{\textup{rank}\:\:}

\newcommand{\GD}{\mathscr{D}}

\newcommand{\present}{\mathscr{P}}

\allowdisplaybreaks[1]

\title{Maximal subgroups of free idempotent generated semigroups over the full transformation monoid}

\author{R. Gray\thanks{Supported by an an EPSRC Postdoctoral Fellowship, and partially supported by FCT and FEDER, project POCTI-ISFL-1-143 of Centro de \'{A}lgebra da Universidade de Lisboa, and by the project PTDC/MAT/69514/2006.},
N. Ru\v{s}kuc\thanks{The second author gratefully acknowledges the financial support and hospitality of 
the Centro de Algebra da Universidade de Lisboa (CAUL) during two visits to 
Lisbon in 2010, financed by FCT and FEDER, during which this research was 
carried out.}}

\begin{document}

\maketitle

\begin{abstract}
Let $T_n$ be the full transformation semigroup of all mappings from the set $\{1, \dots, n\}$ to itself under composition. Let $E = E(T_n)$ denote the set of idempotents of $T_n$ and let $e \in E$ be an arbitrary idempotent satisfying $|\im(e)|=r \leq n-2$. We prove that the maximal subgroup of the free idempotent generated semigroup over $E$ containing $e$ is isomorphic to the symmetric group $S_r$. 

\bigskip

\end{abstract}

\section{Introduction and statement of the Main Theorem} 
\label{sec_intro}

Let $S$ be a semigroup with the set of idempotents $E(S)$ and let $\langle E(S) \rangle$ be the subsemigroup of $S$ generated by $E(S)$. The semigroup $S$ is said to be \emph{idempotent generated} if $S = \langle E(S) \rangle$. 
Idempotent generated semigroups are of interest for a variety of reasons, and consequently they have received considerable attention in the literature. Firstly, they are `generic' in the sense that they satisfy the following universal property: every semigroup can be embedded in an idempotent generated semigroup (see \cite{howie66}). Secondly, many semigroups that occur `in nature' have the property that they are idempotent generated. For example, Howie \cite{howie66} showed that the semigroup of all non-invertible transformations from a finite set to itself is idempotent generated. Erd\"{o}s \cite{erdos67} proved the analogous result for the full linear monoid of all $n \times n$ matrices over an arbitrary field, showing that every non-invertible matrix is expressible as a product of idempotent matrices. This result was later shown to hold more generally for the semigroup of $n \times n$ matrices over an arbitrary division ring; see \cite{laffey83}. Idempotent generated semigroups (and the related theory of biordered sets) also play an important role in the theory of reductive algebraic monoids \cite{putcha88,renner2005}. Relating to this, generalising the abovementioned result of Erd\"{o}s, recently Putcha \cite{putcha06} has obtained necessary and sufficient conditions for a reductive linear algebraic monoid to have the property that every non-invertible element is a product of idempotents.

The \emph{free idempotent generated semigroup over $E$} is the universal object in the category of all 
idempotent generated semigroups whose sets of idemotents
are isomorphic to $E$.
Here $E$ is viewed as a \emph{biordered set}; see \cite{nambooripad79,easdown85,higgins92}.
In fact, in this paper we will not require any theory of  biordered sets.
We just need the definition of the free idempotent generated semigroup over $E$, which is provided by the following presentation:
\begin{equation}
\label{eq1}
\IG(E)=\langle E\:|\: e\cdot f= ef\ (e,f\in E,\ \{ e,f\}\cap\{ef,fe\}\neq\emptyset)\rangle. 
\end{equation}
(It is an easy exercise to show that if, say, $fe\in\{e,f\}$ then $ef\in E$.
In the defining relation $e\cdot f=ef$ the left hand side is a word of length $2$, and $ef$ is the product of $e$ and $f$ in $S$, i.e. a word of length $1$.) 
The idempotents of $S$ and $\IG(E)$ are in the natural $1-1$ correspondence, and we will identify the two sets throughout.
For an idempotent $e\in E$, the maximal subgroup of $\IG(E)$ containing $e$ will be denoted by $H_e$.

Free idempotent generated semigroups were introduced by Nambooripad \cite{nambooripad79}, and the nature of their maximal subgroups
quickly emerged as a key question. In the first phase of the development of the subject several sets of conditions were found which imply freeness of maximal subgroups \cite{mcelwee02,nambooripad80,pastijn77,pastijn80}.
The first example of a non-free maximal subgroup was constructed by Brittenham, Margolis and Meakin in \cite{margolismeakin09},
and was followed by the present authors' construction \cite{grayta}  showing that in fact all groups arise as maximal subgroups of
free idempotent generated semigroups.

In contrast, the structure of the maximal subgroups of of free idempotent generated semigroups arising from naturally occurring biordered sets is still far from being well understood.
In their recent article \cite{margolismeakinip} Brittenham, Margolis and Meakin prove that if $S$ is the full $n\times n$ matrix monoid over the division ring $Q$ and $e\in S$ is an idempotent of rank $1$, then the maximal subgroup $H$ of $\IG(E(S))$ containing $e$ is isomorphic to the multiplicative group of $Q$.
They also conjecture (in Section 5) that the maximal subgroup of an idempotent of rank $r$ with $r\leq n/2$ is isomorphic
to the $r$-dimensional general linear group over $Q$.

The purpose of the present paper is to give a complete description of maximal
subgroups
of the free idempotent generated semigroups arising from finite full transformation semigroups. More precisely, we prove:

\begin{thm}
Let $T_n$ be the full transformation semigroup, let $E$ be its set of idempotents, and let $e \in E$ be an arbitrary idempotent
with image size $r$ ($1\leq r\leq n-2$). Then the maximal subgroup $H_e$ of the free idempotent generated semigroup $\IG(E)$ containing $e$ is isomorphic to the symmetric group $S_r$. 
\end{thm}

If $e$ is the identity mapping then $H_e$ is obviously trivial.
If $|\im e|=n-1$ then $H_e$ is known to be free; this is most easily seen from the presentation \eqref{eqn_top}--\eqref{eqn_bottom} for $H_e$ below, upon observing that trivially there are no singular squares in this case.

The paper is entirely devoted to proving the above result, and is structured as follows: The following section introduces concepts and notation on which the proof of the above theorem will rely, the proof is outlined in Section \ref{sec3}, and executed in Sections \ref{sec4}--\ref{sec8}.
It uses, at least implicitely, methodology introduced in \cite{grayta}.
The success of this approach was somewhat surprising to the authors, as the methods of \cite{grayta} were designed to deal with very different, artificially constructed examples, 
rather than a class of naturally arising examples as is the case here.
Another curious difference is that in \cite{grayta} the connection between the maximal subgroups of $S$ and those of $\IG(E)$ was
in a sense as loose as possible: the former were always trivial, and the latter arbitrary.
This time we have the other extreme: the maximal subgroups containing $e$ in $T_n$ and in $\IG(E)$ are identical!

On one level, the proof of the Main Theorem is simple: 
We take a known general presentation for $H_e$
and use Tietze transformations to turn it into the classical Coxeter presentation for $S_r$.
Technically, however, the argument is considerably more involved. 
This is due to the intricate way in which the Coxeter presentation is encoded in the structure of $T_n$,
reflected primarily in the combinatorics of kernels (partitions) and images (subsets),
and we invite the reader to keep an eye on this, rather beautiful to our mind, aspect of the Main Theorem. 
This encoding is especially subtle when $r$ is large (e.g. equal to $n-2$): Here one can see examples where a generator that needs to be eliminated can be eliminated by one and only one relation, or a relation that needs to be exhibited can
be `read off' from a single kernel/image configuration. This so much so that the authors had several `doubting moments' while working on this project, when they thought that perhaps, as $n$ becomes large, a free factor might creep in! 
On the negative side, the `tightness' of the encoding meant that we have not been able to `distil'
a general criterion for why sometimes the maximal subgroups of a free idempotent generated semigroup happen to be isomorphic to those of the original semigroup. 

\section{Ingredients of the proof} 
\label{sec1}

A presentation for a general maximal subgroup in a general free idempotent generated semigroup was given in \cite[Theorem 5]{grayta}.
Much of what we do below amounts to setting up this general presentation in the particular case under consideration here. This mainly consists in translating from the language of Green's relations of abstract semigroup theory to
that of kernels and images in $T_n$  (see \cite[Exercise 2.6.16]{howie95}) and making concrete choices for entities whose mere existence is asserted in \cite{grayta}.

\subsection{Mappings}
\label{subsec1}

Since this paper is entirely devoted to proving our Main Theorem, the positive integers $n$ and $r$ (with $r\leq n-2$)
and the idempotent transformation $e$ of rank $r$ will be fixed throughout.
Transformations will be written to the right of their arguments and will be composed from left to right.
By \cite[(IG3)]{grayta}, without loss of generality we may assume that
\[
e = 
\begin{pmatrix}
1 & 2 & \dots & r-1 & r & r+1 & \dots & n \\ 
1 & 2 & \dots & r-1 & r & r   & \dots & r  
\end{pmatrix}.
\]
Also we let
\[
D(n,r) = \{ \alpha \in T_n : \rank(\alpha) = r \}.
\]
This, of course, is the $\mathcal{D}$-class of $e$, a fact that is relevant for the subsequent considerations, but one that will remain in the background, as we will be couching our arguments in terms of combinatorics on $T_n$ (sets and partitions), rather than abstract Green's relations.

Let $I$ denote the set of all partitions of $[1,n]$ into $r$ non-empty classes, and let $J$ denote the set of all $r$-element subsets of $[1,n]$. So $I$ and $J$ index the sets of $\mathcal{R}$- and $\mathcal{L}$-classes, respectively, of the $\mathcal{D}$-class $D(n,r)$. 

The $\mathcal{H}$-classes are indexed by the set $I\times J$.
 It is well known (see for instance \cite{howie95}) and very easy to prove that the $\mathcal{H}$-class 
consisting of mappings with kernel $P\in I$ and image $A\in J$ contains an idempotent if and only if $A$ is a transversal of $P$. This will be denoted $A \perp P$, and  $e_{P,A}$ will stand for the unique idempotent in  this $\mathcal{H}$-class.

Presentation \eqref{eq1} introduces a slight notational confusion, which we will maintain throughout the paper:
An element from $E$ can be considered as an idempotent mapping in $T_n$, or as an abstract generator
from the presentation.
Likewise, a word from $E^\ast$ can be considered as a mapping from $T_n$, or as an element of $\IG(E)$.

Let $\lexless$ be the \emph{lexicographic ordering} on the set $J$:
For two sets $A,B\in J$ we write $A\lexless B$ if $A=\{ a_1,\dots, a_r\}$,
$B=\{ b_1,\dots,b_r\}$, with $a_1<\dots <a_r$, $b_1<\dots <b_r$ and, for some $k\in\{1,\dots,r\}$ we have
$a_i=b_i$ ($i=1,\dots,k-1$) and $a_k<b_k$.
This turns $J$ into a linearly ordered set.
For every partition $P=\{P_1,P_2,\dots,P_r\}$ from $I$ we let
$$
A(P)=\{ \min P_1,\dots, \min P_r\},
$$
The following is obvious:

\begin{lemma}
\label{APsmallest}
For any $P\in I$ the set $A(P)$ is the $\lexless$-smallest set $A\in J$ such that $A\perp P$.
\end{lemma}

\subsection{Running example}
\label{running}

We will accompany our theoretical exposition with the concrete example of an idempotent of rank $4$ in the full transformation semigroup on $7$ points.
Thus, here we have $n=7$, $r=4$ and
$$
e=\begin{pmatrix}
1 & 2 & 3&4&5&6&7 \\ 
1 & 2 & 3&4&4&4&4  
\end{pmatrix}.
$$
The index set $I$ has $S(7,4)=350$ elements (the Stirling number of the second kind),
while  $J$ has $\binom{7}{4}=35$ elements.
There are $2240$ pairs $(P,A)$ with $A\perp P$.

\subsection{Permutations}
\label{perms}

All permutations that we encounter in this paper will come from the symmetric group $S_r$.
We will use both  cycle notation and  image notation for permutations.
The former will be indicated by using round braces, and the latter by using square braces. 
Thus, for example, the permutation 
$$
\pi =\begin{pmatrix}
1&2&3&4 \\ 
3&2&4&1  
\end{pmatrix}
$$
may be written as $(1\: 3\: 4)$ or as $[3,2,4,1]$.

We will also make use of the following measure of complexity of permutations related to the weak 
Bruhat order
 (see \cite[Chapter 2]{bjorner}).
Suppose that $\pi=[p_1,\dots,p_r]\in S_r$.
An entry $p_k$ is said to be the \emph{descent start} in $\pi$ if there exists $l>k$ such that
$p_k>p_l$; any such pair $(p_k,p_l)$ is called a \emph{descent}.
The number of descent starts in $\pi$ will be called its \emph{descent number}, and will be denoted by $\Delta(\pi)$.
For example, $\Delta(3241)=3$, since $3$, $2$ and $4$ are all descent starts.

\begin{sloppypar}
We record the following obvious descriptions of permutations with small descent numbers:
\end{sloppypar}

\begin{lemma}
\label{smalldescent}
For $\pi\in S_r$ the following hold:
\begin{enumerate}[label=\textup{(\arabic*)}, leftmargin=*,widest=2]
\item
$\Delta(\pi)=0$ if and only if $\pi=()$.
\item
$\Delta(\pi)=1$ if and only if $\pi$ is a contiguous cycle of the form
\[
\pi=(k+l\ k+l-1\ \dots\ k),
\]
for some $k$ and $l$ with $1\leq k < k+l \leq r$.
\end{enumerate}
\end{lemma}

\subsection{Actions and Schreier representatives}
\label{subsec2}

There is a natural action of $T_n$ on the set of all $\mathcal{H}$-classes in the $\mathcal{R}$-class of $e$ with $0$ adjoined.
This action is naturally equivalent to the action of $T_n$ on $J \cup \{ 0 \}$ given by:
\[
A \cdot \alpha = 
\begin{cases}
A \alpha = \{ a \alpha : a \in A \} ,
& \mbox{if} \ |A\alpha|=r ,\\
0,
& \mbox{if} \ |A\alpha|<r.  
\end{cases}
\]  

A set of words $\rho_A, \rho_A' \in E^*$ form a \emph{Schreier system of representatives} if the following conditions are satisfied:
\begin{enumerate}[label=\textsf{(S\arabic*)}, leftmargin=*,widest=2]
\item
$\rho_A$ and $\rho_A'$ induce mutually inverse bijections between $[1,r]$ and $A$;
\item
every prefix of $\rho_A$ (including the empty word) is equal to some $\rho_B$. 
\end{enumerate}

We define a particular Schreier system inductively with respect to $\lexless$ as follows.
First we set:
\[
\rho_{[1,r]} = \rho_{[1,r]}^{-1} = \epsilon,
\]  
the empty word. 
Now let $A=\{a_1,\dots,a_r\}\neq [1,r]$, with $a_1<\dots<a_r$,
and let $m\in [1,r]$ be the smallest subscript such that $a_m\neq m$; note that $a_m>1$ and $a_m-1\not\in A$.
Let $B\in J$ and $P=\{ P_1,\dots, P_r\}\in I$ be given by:
\begin{eqnarray*}
&&B=\{ a_1,\dots, a_{m-1}, a_m-1,a_{m+1},\dots, a_r\},\\
&&P=\{ [1,a_1],[a_1+1,a_2],\dots, [a_{r-2}+1,a_{r-1}],[a_{r-1}+1,n]\},
\end{eqnarray*}
and define
\[
\rho_A=\rho_B e_{P,A},\ \rho_A^\prime=e_{P,B} \rho_B^\prime.
\]
This definition makes sense since $B\lexless A$ and 
$\{A,B\}\perp P$
(meaning both $A$ and $B$ are transversals for $P$).
We record, for future use, that when considered as a mapping in $T_n$,  the word $\rho_A$ restricted to $[1,r]$ is the unique order preserving bijection from $[1,r]$ to $A$; 
this follows by an easy inductive argument.

\subsection{Labels}
\label{label} %:-)

Suppose that $P = \{ P_1, \dots, P_r   \} \in I$ and $A\in J$ are such that $A\perp P$.
Recall that also $A(P)\perp P$. Let $\gamma_{P,A}: A(P)\rightarrow A$ be the bijection
which sends $\min P_i$ to the unique element of $A$ belonging to $P_i$.
Recall that $\rho_{A(P)}\rrestriction_{[1,r]}$ and $\rho_A\rrestriction_{[1,r]}$ are the unique order preserving bijections from
$[1,r]$ onto $A(P)$ and $A$ respectively.
Define the \emph{label} of $(P,A)$ to be
\begin{equation}
\label{labeleq}
\lambda(P,A) = \rho_{A(P)}\rrestriction_{[1,r]} \gamma_{P,A} \rho_A^{-1},
\end{equation}
which is clearly a permutation from $S_r$.

An immediate corollary of this definition is the following:

\begin{lemma}
\label{firstlabel}
For every $P\in I$ we have $\lambda(P,A(P))=()$.
\end{lemma}

\begin{runex}
Let us compute the label $\lambda(P,A)$ where 
\[
P=\{ \{1\},\{2,3,5\},\{4,7\},\{6\}\},\ A=\{1,4,5,6\}.
\]
Clearly $A(P)=\{1,2,4,6\}$, and so
\[
\rho_{A(P)}\rrestriction_{[1,4]}=\begin{pmatrix} 1&2&3&4\\ 1&2&4&6\end{pmatrix},\ 
\gamma_{P,A}=\begin{pmatrix} 1&2&4&6\\ 1&5&4&6\end{pmatrix},\ 
\rho_A^{-1}\rrestriction_A=\begin{pmatrix} 1&4&5&6\\ 1&2&3&4\end{pmatrix}.
\]
The label is
\[
\lambda(P,A)=\rho_{A(P)}\rrestriction_{[1,4]} \gamma_{P,A}\rho_A^{-1}
=\begin{pmatrix} 1&2&3&4\\ 1&3&2&4\end{pmatrix}= (2\ 3).
\]
\end{runex}

\begin{remark}
In practice one can compute the label of an arbitrary pair $A \perp P$ quickly and easily as
follows. First write:
\begin{eqnarray*}
&&A=\{a_1,\dots,a_r\}, \ a_1<\dots<a_r,\\
&&  P=\{P_1,\dots,P_r\}, \ \min P_1<  \dots<  \min P_r.
\end{eqnarray*}
Then write out the sets $P_1$ up to $P_r$ in order, and underneath each $P_i$
write the unique element $a_{l_i}$ from $A$ that belongs to $P_i$, giving:
\[
\begin{pmatrix}
P_1&  P_2&  \ldots&  P_r \\
a_{l_1}&  a_{l_2}&  \ldots&  a_{l_r}
\end{pmatrix}   .
\]
Then the label is given by
keeping the subscripts:
\[
\lambda(P,A) =
\begin{pmatrix}
1&  2&  \ldots&  r \\
 l_1&  l_2&  \ldots&  l_r
\end{pmatrix}.
\]
\end{remark}

\begin{remark}
For another viewpoint, it may be checked that the labels correspond to the 
non-zero entries of the structure matrix of a certain natural Rees matrix 
representation of the principal factor of $T_n$ arising from the $\GD$-class of $e$, as 
described in \cite[Section 3.2]{howie95}.
\end{remark}

\subsection{Singular squares}
\label{subsec3}

Modifying the notation (but not the substance) from \cite[Section 3]{grayta},
we call a quadruple $(P,Q,A,B)\in I\times I\times J\times J$ a \emph{square} if $\{A,B\} \perp \{P,Q\}$. A square is \emph{singular} if, in addition, there exists $e \in E$ such that either:
\begin{align}
& e e_{P, A} = e_{P, A}, \quad
e e_{Q, A} = e_{Q, A}, \quad
e_{P,A} e = e_{P,B}, \quad
e_{Q,A} e = e_{Q,B}, \; or \label{eqnLR} 
\\
& e_{P,A}e = e_{P,A}, \quad
e_{P,B}e = e_{P,B}, \quad
e e_{P,A} = e_{Q,A}, \quad
e e_{P,B} = e_{Q,B}. \label{eqnUD}
\end{align}
Let $\Sigma_{LR}$ (respectively $\Sigma_{UD}$) be the set of all singular squares for 
which condition \eqref{eqnLR}
(resp. \eqref{eqnUD}) holds, and
let $\Sigma=\Sigma_{LR}\cup\Sigma_{UD}$, the set of all singular squares.
We call the members of $\Sigma_{LR}$ the \emph{left-right} singular squares, and those of $\Sigma_{UD}$ the \emph{up-down} singular squares.

\begin{lemma}
\label{singularity}
The following conditions are equivalent for a square $(P,Q,A,B)$:
\begin{enumerate}[label=\textup{\textsf{(SQ\arabic*)}}, leftmargin=*,widest=3]
\item
\label{SQ1}
$(P,Q,A,B)$ is singular;
\item
\label{SQ2}
$
\{ (A \cap P_i, B \cap P_i): i=1, \dots, r \}
=
\{ (A \cap Q_i, B \cap Q_i): i=1, \dots, r \};
$
\item
\label{SQ3}
$
\lambda(P,A)^{-1}
\lambda(P,B)
=
\lambda(Q,A)^{-1}
\lambda(Q,B).
$ 
\end{enumerate}
\end{lemma}

\begin{proof}
\ref{SQ1}$\Rightarrow$\ref{SQ2}
Suppose $a\in A\cap P_i$, $b\in B\cap P_i$,
$a\in Q_j$, $b\in Q_k$.
We aim to show that $j=k$.
If $(P,Q,A,B)\in\Sigma_{LR}$
with an idempotent $e$ satisfying \eqref{eqnLR}, then from
$$
ae=ae_{P,A}e=ae_{P,B}=b
$$
we have
$$
Q_j\ni a=ae_{Q,A}=aee_{Q,A}=be_{Q,A}\in Q_k,
$$
implying $j=k$. Similarly, if $(P,Q,A,B)\in\Sigma_{UD}$, with
$e$ satisfying \eqref{eqnUD}, then from
$$
ae=ae_{P,A}e=ae_{P,A}=a
$$
we have
$$
Q_k\ni b=ae_{P,B}=aee_{P,B}=ae_{Q,B}\in Q_j,
$$
and so again $j=k$.
\medskip

\ref{SQ2}$\Rightarrow$\ref{SQ1}
Define a mapping $e\in T_n$ by
$$
xe=
\begin{cases}
y & \text{if } B\cap P_i=\{x\},\ A 
\cap Q_i=\{y\} \text{ for some } i,\\
x & \text{otherwise.}
\end{cases}
$$
A routine verification shows that $e$ is an idempotent and that \eqref{eqnLR} is satisfied.
Thus $(P,Q,A,B)$ is a (left-right) singular square.
\medskip

\ref{SQ2}$\Leftrightarrow$\ref{SQ3}
Using \eqref{labeleq}  we have
\begin{eqnarray}
\nonumber
&&\lambda(P,A)^{-1}
\lambda(P,B)
=
\lambda(Q,A)^{-1}
\lambda(Q,B)
\\
\nonumber
&\Leftrightarrow& (\rho_{A(P)}\rrestriction_{[1,r]}\gamma_{P,A} \rho_A^{-1})^{-1}
\rho_{A(P)}\rrestriction_{[1,r]} \gamma_{P,B} \rho_B^{-1}\\
\nonumber
&&
\hspace{25mm}=
(\rho_{A(Q)}\rrestriction_{[1,r]} \gamma_{Q,A} \rho_A^{-1})^{-1}
\rho_{A(Q)} \rrestriction_{[1,r]} \gamma_{Q,B} \rho_B^{-1}
\\
\label{cond1}
&\Leftrightarrow&
\gamma_{P,A}^{-1}\rrestriction_A \gamma_{P,B} = 
\gamma_{Q,A}^{-1}\rrestriction_A \gamma_{Q,B}.
\end{eqnarray}
Note that the mapping $\gamma_{P,A}^{-1}\rrestriction_A \gamma_{P,B}$
maps the only element of $A\cap P_i$, via $\min P_i$, to the only element of $B\cap P_i$.
An analogous statement holds for $\gamma_{Q,A}^{-1}\rrestriction_A \gamma_{Q,B}$.
Hence \eqref{cond1} is equivalent to \ref{SQ2}, as required.
\end{proof}

\begin{remark}
It follows from the above proof that every singular square in $T_n$ is an LR-square.
\end{remark}

\begin{remark}
Suppose a square $(P,Q,A,B)$ is a \emph{rectangular band}, meaning that the set
$\{ e_{P,A},e_{P,B},e_{Q,A},e_{Q,B}\}$ is closed under multiplication.
It is an easy exercise to show that this occurs if and only if $e_{P,A}e_{Q,B}=e_{P,B}$.
Suppose that we have $a\in A$ and $b\in B$ such that $a,b\in P_i$ for some $i\in [1,r]$.
Further, suppose that $a\in Q_j$, while $b\in Q_k$.
Since $a\in A$ it follows that $ae_{P,A}=a$,
while from $a,b\in P_i$ and $b\in B=\im e_{P,B}$ it follows that $a e_{P,B}=b$.
Hence $ae_{Q,B}=ae_{P,A}e_{Q,B}=b$, and, since $e_{Q,B}$ preserves the blocks of $Q$,
it follows that $j=k$.
This demonstrates that the condition \ref{SQ2} of Lemma \ref{singularity} is satisfied,
proving that every rectangular band is singular.
This is the analogue for $\IG(E(T_n))$ of \cite[Theorem 4.3]{margolismeakinip} which proves the same fact for the free
idempotent generated semigroup over a full matrix monoid over a division ring.
\end{remark}

\begin{runex}
Let
\begin{align*}
& P=\{ \{1\},\{2,3,5\},\{4,7\},\{6\}\}, && Q=\{\{1\},\{2,3,6\},\{4,7\},\{5\}\},\\
 & A=\{1,4,5,6\}, && B=\{1,5,6,7\}.
\end{align*}
We saw in Subsection \ref{label} that $\lambda(P,A)=(2\ 3)$. Performing the same calculations for the other three pairs, we obtain
\[
\lambda(P,B)=(3\ 4),\ \lambda(Q,A)=(2\ 4\ 3),\ \lambda(Q,B)=(2\ 3\ 4).
\]
Since 
\[
(2\ 3)^{-1} (3\ 4)=(2\ 4\ 3)=(2\ 4\ 3)^{-1} (2\ 3\ 4),
\]
and so the square $(P,Q,A,B)$ is singular.
By way of contrast, for
\begin{align*}
& P^\prime=\{ \{1\},\{2,4\},\{3,6\},\{5,7\}\}, && Q^\prime=\{\{1\},\{2,6,7\},\{3,5\},\{4\}\},\\
 & A^\prime=\{1,3,4,7\}, && B^\prime=\{1,4,5,6\},
\end{align*}
we have
\[
\lambda(P^\prime,A^\prime)= (2\ 3),\ \lambda(P^\prime,B^\prime)=(3\ 4),\ \lambda(Q^\prime,A^\prime)=(2\ 4\ 3),\ \lambda(Q^\prime,B^\prime)=(2\ 4).
\]
Since
\[
(2\ 3)^{-1} (3\ 4)=(2\ 4\ 3)\neq (2\ 3)=(2\ 4\ 3)^{-1}(2\ 4),
\]
the square $(P^\prime,Q^\prime,A^\prime,B^\prime)$ is not singular.
\end{runex}

\begin{remark}
It is not true that all permutations from $S_r$ arise as labels of singular squares, although this \emph{is} the case
in our running example ($n=7$, $r=4$). For example, computation using \textsf{GAP} (\cite{GAP4}) shows that for $n=7$, $r=5$
only $46$ out of $120$ permutations are labels of singular squares.
It is, however, true that all Coxeter transpositions are always present as labels. This follows from our argument below, but is also fairly straightforward to prove by a direct construction.
\end{remark}

\section{Outline of the proof}
\label{sec3}

The collection of all singular squares yields a presentation for $H$.
This was first proved by Nambooripad \cite{nambooripad79} for the regular case,
and extended in \cite{grayta} to the non-regular case. 
Following
\cite{grayta}, we have that $H$ is defined by the presentation with generators
$$
F=\{ f_{P,A} \st P\in I,\ A\in J,\ A\perp P\},
$$
and the defining relations
\begin{alignat}{2}
\label{eqn_top}
&
f_{P,A}=f_{P,B} &\quad & (A \perp P, \; B \perp P, \; \rho_A e_{P,B} = \rho_{A \cdot e_{P,B}}),
\\
\label{eqn_middle}
& f_{P,A(P)}=1 && (P\in I),
\\
\label{eqn_bottom}
& f_{P,A}^{-1}f_{P,B}=f_{Q,A}^{-1}f_{Q,B} && ((P,Q,A,B)\in\Sigma).
\end{alignat}
Let us denote this presentation by $\present$. The generator $f_{P,A}$ represents the element
$e \rho_{A(P)} e_{P,A} \rho_A^{-1}$ in $\IG(E)$.
If this element is interpreted as a transformation in $T_n$, it belongs to the copy of the symmetric group $S_r$
consisting of all mappings with kernel $P$ and image $A$.
Interpreted further as an element of the natural copy of  $S_r$ (acting on $[1,r]$), 
via the identification $a_i\mapsto i$, where $A=\{ a_1,\dots,a_r\}$, $a_1<\dots<a_r$,
this element is equal to the label $\lambda(P,A)$. 
Motivated by this, we extend the scope of the labelling function $\lambda$ to generators from $F$
by setting $\lambda(f_{P,A})=\lambda(P,A)$.

The proof of our main theorem will consist in applying Tietze transformations to the above presentation to eventually obtain the well known \emph{Coxeter presentation}
\begin{alignat}{2}
\nonumber
\langle g_1,\dots, g_{r-1}  \: |\: & 
g_i^2=1 &\quad &  (i\in [1,r-1]),
\\
\label{Cox}
& g_i g_j = g_j g_i && (i,j\in [1,r-1],\ |i-j|>1),
\\
\nonumber
&g_i g_{i+1} g_i = g_{i+1} g_i g_{i+1} && (i\in [1,r-1])\rangle,
\end{alignat}
which defines $S_r$ in terms of the generating set consisting of \emph{Coxeter transpositions}
$(i\ i+1)$ ($i\in [1,r-1]$).
This will be organised as follows:

\begin{enumerate}[label=\textsf{(P\arabic*)}, leftmargin=*,widest=2]
\item
\label{PF1}
 We begin by showing how to eliminate all generators $f_{P,A}$ with the identity label, by showing that 
if $\lambda(f_{P,A})=()$ then $f_{P,A}=1$ is a consequence of
the presentation $\present$.
(Section \ref{sec4}.)
\item
\label{PF2}
If two generators $f_{P,A}$ and $f_{Q,B}$ belong to the same row or column (i.e. $P=Q$ or $A=B$)
and if their labels are equal, then $f_{P,A}=f_{Q,B}$ is a consequence of $\present$.
These are auxiliary results, which are used in the subsequent inductive arguments.
(Section \ref{sec5}.)
\item
\label{PF3}
Next we deal with  the generators whose labels have descent number 1 (reverse contiguous cycles, see Lemma \ref{smalldescent}).
This of course includes those labelled by Coxeter transpositions.
We show that any two such generators with equal labels are themselves equal as a consequence of $\present$
and also how to eliminate those whose label has length greater than $2$ (i.e. it is not a Coxeter transposition).
(Section \ref{sec6}.)
\item
\label{PF4}
Then we eliminate all generators whose labels have descent number greater than $2$ by expressing them as products of generators with smaller descent numbers. (Section \ref{sec7}.)
\item
\label{PF5}
At this stage we are left with a presentation with generators in one-one correspondence with the Coxeter generators of $S_r$,
and which defines a homomorphic pre-image of $S_r$. Thus the following step completes the proof.
\item
\label{PF6}
All the Coxeter relations are consequences of $\present$. (Section \ref{sec8}.)
\end{enumerate}

Our methodology is strongly influenced by that of our previous paper \cite{grayta}. 
In particular, implicit in our arguments below is the use of certain special types of singular squares introduced in \cite{grayta}. 
Modifying slightly the terminology from \cite{grayta}, we say that a
singular square $(P,Q,A,B)$ is: 
\begin{itemize}
\item
a \emph{corner} square if we have already proved that three of the associated generators $f_{P,A}$, $f_{P,B}$, $f_{Q,A}$, $f_{Q,B}$
are equal to $1$, in which case we may deduce that the remaining generator also equals $1$;
\item
a \emph{flush left} (resp. \emph{top}) square if we have already proved that $f_{P,A}=f_{Q,A}$ (resp. $f_{P,A}=f_{P,B}$), 
in which case we may deduce that the remaining two generators are equal;
\item
a 3/4-\emph{square} if $f_{P,A}$ has been shown to equal $1$, in which case we may deduce that $f_{Q,A}f_{P,B}=f_{Q,B}$.
\end{itemize}

In this terminology, the proof of \ref{PF1} can be interpreted as showing that we can start from generators involved in  the relations 
\eqref{eqn_top} and \eqref{eqn_middle} of $\present$, and reach every generator labelled by $()$ via a sequence of corner squares.
To prove \ref{PF2} we show that any pair of generators labelled by the same contiguous cycle can be linked by a sequence of flush squares
while the eliminations under \ref{PF3} and \ref{PF4} are achieved by means of 3/4-squares.
By way of contrast, to prove \ref{PF6} we resort to more complicated types of squares, or, indeed, combinations of overlapping squares.

Since the defining relations \eqref{eqn_bottom} are labelled by singular squares, each time we want to make use of such a relation we need to demonstrate singularity of a square. This is done by computing the relevant labels and checking condition \ref{SQ3}
of Lemma \ref{singularity}.
Quite a few such computations are necessary in our argument, and, 
since they  all follow the same routine pattern, 
we have omitted them except for a few sample ones. 
Probably the most instructive of these is performed in the proof of Lemma \ref{smallestA}.

\begin{sloppypar}
We received the first intimations of the truth of 
the Main Theorem
through `experimental' investigations using \textsf{GAP}.
Of particular help was the amazingly functional Tietze Transformations routine, which could handle huge presentations
arising from \eqref{eqn_top}--\eqref{eqn_bottom} and transform them into `human friendly' ones.
The Tietze routine is a part of the main \textsf{GAP} distribution, and is in the manual creditted back to the work of Havas, Robertson et al. \cite{havas69,havas84,robertson88}.
\end{sloppypar}

\section{Generators with  label \boldmath{$()$}} 
\label{sec4}

Our first step in the proof of the Main Theorem is to eliminate generators with label $()$. We do this in three steps.

\begin{lemma}
\label{fixlemma1a}
Let $a\in [r,n]$, 
$P=\{ \{1\},\dots,\{r-1\},[r,n]\}$, and \\
$A=\{1,\dots,r-1,a\}$.
Then the relation $f_{P,A}=1$ is a consequence of $\present$.
\end{lemma}

\begin{proof}
For $a=r$ we have $A=A(P)$ and $f_{P,A}=1$ belongs to \eqref{eqn_middle}.
Proceeding inductively, suppose that $a\in [r+1,n]$,  and that $f_{P,B}=1$ for $B=\{ 1,\dots,r-1,a-1\}$.
By Subsection \ref{subsec2} we have $\rho_A=\rho_B e_{P,A}$, and hence 
the
relation $f_{P,A}=f_{P,B}$ is in \eqref{eqn_top}.
\end{proof}

We say that a partition $P\in I$ is \emph{convex} if all its classes are intervals.
Clearly, if $P$ is convex and $A\perp P$ we must have $\lambda(P,A)=()$.

\begin{lemma}
\label{fixlemma1b}
If $P$ is a convex partition and $A\perp P$ then the relation $f_{P,A}=1$ is a consequence of $\present$.
\end{lemma}

\begin{proof}
We prove the lemma by a double induction.
The first induction is on $A(P)$ with respect to $\lexless$,
the anchor being provided by Lemma \ref{fixlemma1a}.
Let
\[
P=\{ [p_i,p_{i+1}-1]\st i\in [1,r]\},
\]
where $1=p_1<p_2<\dots<p_r<p_{r+1}=n+1$, and
\[
A(P)=\{ p_1,\dots,p_r\}\neq [1,r].
\]
Let $m\in [1,r-1]$ be the smallest subscript such that $p_{m+1}\neq m+1$.

The second induction is on $A$, again with respect to $\lexless$.
The anchor here is $A=A(P)$, in which case $f_{P,A}=1$ belongs to \eqref{eqn_middle}.
Consider now
\[
A=\{ a_1,\dots,a_r\} \neq A(P), \text{ with } a_1<\dots<a_r.
\]
Note that $A\perp P$ means that $a_i\in [p_i,p_{i+1}-1]$ for all $i\in [1,r]$.
Let $t\in [1,r]$ be the smallest subscript such that $a_t\neq p_t$.
Clearly, $t\geq m$, and we distinguish two cases.

\medskip\noindent
\textit{Case 1: $t=m$.}
We have
\[
a_1=p_1=1,\dots, a_{m-1}=p_{m-1}=m-1, a_m> p_m=m.
\]
Let
\[
\begin{split}
&B=\{ 1,\dots, m-1,a_m-1,a_{m+1},\dots, a_r\},\\
&Q=\{ \{1\},\dots,\{m-1\}, [m,a_m],[a_m+1,a_{m+1}],\\
&\quad\quad\quad\quad\quad\quad\quad\quad\quad\quad\quad\quad\quad
\dots, [a_{r-2}+1,a_{r-1}],[a_{r-1}+1,n]\}.
\end{split}
\]
It is easy to see that $\{A,B\}\perp\{P,Q\}$.
Since both $P$ and $Q$ are convex, all four labels equal $()$, and so the square $(P,Q,A,B)$ is singular by Lemma \ref{singularity}.
Thus we have the relation 
\begin{equation}
\label{y1}
f_{P,A}^{-1}f_{P,B}=f_{Q,A}^{-1}f_{Q,B}
\end{equation}
in \eqref{eqn_bottom}.
The definition of Schreier representatives from Subsection \ref{subsec2} gives $\rho_A=\rho_B e_{Q,A}$, so that the relation $f_{Q,A}=f_{Q,B}$ is in \eqref{eqn_top}.
Also, $B\lexless A$, and by the second induction we have $f_{P,B}=1$. Substituting this into \eqref{y1} yields $f_{P,A}=1$ as required.

\medskip\noindent
\textit{Case 2: $t>m$.}
Let
\begin{multline*}
Q=\{ \{1\},\dots,\{m-1\},[m,p_{m+1}-2],[p_{m+1}-1,p_{m+2}-1],\\
[p_{m+2},p_{m+3}-1],\dots, [p_r,n]\}.
\end{multline*}
Effectively, $Q$ is obtained from $P$ by moving the element $p_{m+1}-1$ from the block $P_m$ into the block $P_{m+1}$.
Thus $A(Q)\lexless A(P)$.
It is straightforward to check that $\{A,A(P)\}\perp \{P,Q\}$.
Since both $P$ and $Q$ are convex, the square $(P,Q,A,A(P))$ is singular, yielding
$f_{P,A}^{-1}f_{P,A(P)}=f_{Q,A}^{-1} f_{Q,A(P)}$.
By the first induction we have $f_{Q,A}=f_{Q,A(P)}=1$, while the relation $f_{P,A(P)}$ is in \eqref{eqn_middle}.
Combining we obtain $f_{P,A}=1$, and the proof is complete.
\end{proof}

\begin{lemma}
\label{fixlemma1}
If $P\in I$, $A\in J$ are such that $A\perp P$ and $\lambda(P,A)=()$ then the relation
$f_{P,A}=1$ is a consequence of $\present$.
\end{lemma}

\begin{proof}
We consider $P$ fixed and induct on $A$ with respect to $\lexless$.
When $A=A(P)$ (the anchor) the relation $f_{P,A}=1$ is in \eqref{eqn_middle}.
Let now
\begin{eqnarray*}
&&A=\{a_1,\dots,a_r\}\neq A(P),\ a_1<\dots<a_r,\\
&& P=\{P_1,\dots,P_r\},\ p_1=\min P_1< \dots< p_r=\min P_r.
\end{eqnarray*}
From $\lambda(P,A)=()$ it follows that $a_i\in P_i$ for all $i\in [1,r]$.
Let $m\in [1,r]$ be the smallest subscript for which $p_m\neq a_m$. 
Define:
\begin{eqnarray*}
B&=&\{ p_1,\dots, p_m,a_{m+1},\dots,a_r\},\\
Q&=&\{ [p_1,p_2-1],\dots, [p_{m-1},p_m-1], [p_m,a_{m+1}-1],[a_{m+1},a_{m+2}-1],\\
&&\hspace{60mm}\dots,
[a_{r-2},a_{r-1}-1],[a_{r-1},n]\}.
\end{eqnarray*}
It is easy to see that $\{A,B\}\perp \{P,Q\}$, and that all four labels equal $()$.
By Lemma \ref{singularity} the square $(P,Q,A,B)$ is singular, and we have the relation $f_{P,A}^{-1} f_{P,B}=f_{Q,A}^{-1} f_{Q,B}$ in \eqref{eqn_bottom}.
Since $Q$ is convex we have $f_{Q,A}=f_{Q,B}=1$ by Lemma \ref{fixlemma1b}.
We also have $f_{P,B}=1$ by induction, since $B\lexless A$.
It follows that $f_{P,A}=1$, as required.
\end{proof}

\section{Generators in the same row or column}
\label{sec5}

The aim in this section is to prove that generators with equal labels which lie in the same row or column are equal as
a consequence of $\present$. We begin with an auxiliary result:

\begin{lemma}
\label{lemmav1}
Let $A\in J$ and $P,Q\in I$ be such that $A\perp P$, $Q\perp A$ and $\lambda(P,A)=\lambda(Q,A)$.
Furthermore, suppose $P=\{P_1,\dots,P_r\}$ with
$\min P_1 <\dots <\min P_r$ and $Q=\{ Q_1,\dots,Q_r\}$ with
$\min Q_1<\dots <\min Q_r$.
If there is an increasing sequence $b_1<\dots<b_r$ with
$b_i\in P_i\cap Q_i$ ($i=1,\dots,r$) then $f_{P,A}=f_{Q,A}$ is a consequence of $\present$.
\end{lemma}

\begin{proof}
Let $B=\{ b_1,\dots, b_r\}$; clearly $B\perp P$, $B\perp Q$, and
$\lambda(P,B)=\lambda(Q,B)=()$.
By Lemma \ref{singularity} the square $(P,Q,B,A)$ is singular, and so we have the relation
$f_{P,B}^{-1}f_{P,A}=f_{Q,B}^{-1}f_{Q,A}$.
By Lemma \ref{fixlemma1} we have $f_{P,B}=f_{Q,B}=1$, leaving us with
$f_{P,A}=f_{Q,A}$ as desired.
\end{proof}

\begin{lemma}
\label{lemmav2}
If $A\in J$ and $P,Q\in I$ are such that $A\perp \{P,Q\}$ and
$\lambda(P,A)=\lambda(Q,A)$ then $f_{P,A}=f_{Q,A}$ is a consequence of $\present$.
\end{lemma}

\begin{proof}
As usual, without loss of generality we may suppose that
\begin{eqnarray*}
&&
A=\{ a_1,\dots,a_r\},\ a_1<\dots<a_r,
\\
&&
P=\{P_1,\dots,P_r\},
\ p_i=\min P_i\ (i=1,\dots,r),\ p_1<\dots<p_r,
\\
&&
Q=\{ Q_1,\dots,Q_r\},
\ q_i=\min Q_i\ (i=1,\dots,r),\ q_1<\dots<q_r,
\\
&&
\lambda(P,A)=\lambda(P,B)=\pi\in S_r,
\\
&&
\{ a_{i\pi} \} =A\cap P_i\cap Q_i\ (i=1,\dots, r).
\end{eqnarray*}
Let $u+1$ be the smallest index in which the sequences $(p_1,\dots,p_r)$ and
$(q_1,\dots,q_r)$ differ;
in other words $p_i=q_i$ ($i=1,\dots,u$) and $p_{u+1}\neq q_{u+1}$.
Clearly $u\geq 1$ as $p_1=q_1=1$.
Define the \emph{distance} $d(P,Q)$ between $P$ and $Q$ to equal $r-u$.
We prove the lemma by induction on $d(P,Q)$.

In the anchor case $d(P,Q)=0$ we have $p_i=q_i$ for all $i=1,\dots,r$, and so $f_{P,A}=f_{Q,A}$ 
follows from Lemma \ref{lemmav1}.
Suppose $d(P,Q)=d=r-u>0$ and that the lemma holds for all pairs of partitions at a smaller distance
from each other.
Suppose without loss of generality that $p_{u+1}< q_{u+1}$.
The element $p_{u+1}$ belongs to some $Q_v$.
Since $p_{u+1}<q_{u+1}$, we must have $v\leq u$.
But then we must have $p_{u+1}\neq p_v=q_v$.
Furthermore, $p_{u+1}\in P_{u+1}\setminus Q_{u+1}$ implies $p_{u+1}\not\in A$.

Transform $Q$ into a new partition $R=\{R_1,\dots, R_r\}$ by moving $p_{u+1}$ from $Q_v$ to $Q_{u+1}$:
\[
R_v=Q_v\setminus \{ p_{u+1}\},\ R_{u+1}=Q_{u+1}\cup\{ p_{u+1}\},\ 
R_i=Q_i\ (i\neq v,u+1).
\]
Since still $a_{i\pi}\in R_i$ for all $i=1,\dots ,r$ it follows that
$\lambda(R,A)=\pi$.
Furthermore, $d(P,R)\leq d-1$ (because the minimum of $R_{u+1}$ is $p_{u+1}$, the same as for $P_{u+1}$), and by induction we have $f_{P,A}=f_{R,A}$.
Finally, the increasing sequence $q_1<\dots<q_r$ satisfies $q_i\in Q_i\cap R_i$ and so
$f_{Q,A}=f_{R,A}$ follows from Lemma \ref{lemmav1}.
Combining we obtain $f_{P,A}=f_{Q,A}$, completing the proof.
\end{proof}

\begin{lemma}
\label{lemmav2a}
If $A,B\in J$ and $P\in I$ are such that $\{A,B\}\perp P$ and
$\lambda(P,A)=\lambda(P,B)$ then $f_{P,A}=f_{P,B}$ is a consequence of $\present$
\end{lemma}

\begin{proof}
Suppose $P=\{P_1,\dots,P_r\}$, $A=\{a_1,\dots,a_r\}$, $B=\{ b_1,\dots,b_r\}$,
with
\[
\min P_1 <\dots < \min P_r,\ 
a_1<\dots<a_r,\ 
b_1<\dots<b_r.
\]
Note that $\lambda(P,A)=\lambda(P,B)$ means that for all $i,j\in [1,r]$ we have
$$
a_i \in P_j \Leftrightarrow b_i\in P_j.
$$
Define a new partition $Q=\{ Q_1,\dots, Q_r\}$ by
\[
Q_i=\{a_i,b_i\} \ (i\in [2,r]),\ 
Q_1=[1,n]\setminus (Q_2\cup\dots\cup Q_r).
\]
Let $q_i=\min Q_i$, so that $A(Q)=\{ q_1,\dots,q_r\}$, and note that $q_1<\dots<q_r$.
It is clear that $\{A,B\}\perp Q$. Furthermore,
we have $\rho_{A(Q)}\st i\mapsto c_i$,
$\gamma_{Q,A}\st c_i\mapsto a_i$, $\gamma_{Q,B}\st c_i \mapsto b_i$,
$\rho_A\st i\mapsto a_i$, $\rho_B\st i\mapsto b_i$,
and it readily follows that $\lambda(Q,A)=\lambda(Q,B)=()$.
By Lemma \ref{singularity}, the square $(P,Q,A,B)$ is singular, and so the relation
$f_{P,A}^{-1}f_{P,B}=f_{Q,A}^{-1}f_{Q,B}$ is in \eqref{eqn_bottom}.
By Lemma \ref{fixlemma1}, the relations $f_{Q,A}=f_{Q,B}=1$ are consequences of $\present$,
leaving us with $f_{P,A}=f_{Q,A}$, as required.
\end{proof}

\begin{remark}
The above proof shows that any pair of generators that belong to the same row and have equal labels belong to a singular square
the other two vertices of which are labelled $()$.
The analogous assertion is not true for pairs of generators that belong to the same column, 
in which case the proof of Lemma \ref{lemmav2} merely asserts
that such a pair can be linked by a sequence of appropriate squares.
\end{remark}

\begin{runex}
Let $P=\{ \{1\},\{2\},\{3,6\},\{4,5,7\}\}$, $A=\{ 1,2,5,6\}$ and $B=\{1,2,4,6\}$. 
The generators $f_{P,A}$ and $f_{P,B}$ are in the same row, and both have label $(3\ 4)$.
If $P^\prime=\{ \{ 1,3,7\},\{2\},\{4,5\},\{6\}\}$ then $\lambda(P^\prime,A)=\lambda(P^\prime,B)=()$.
Now let $Q=\{\{1\},\{2,3\},\{4,6,7\},\{5\}\}$; we have $\lambda(Q,A)=(3\ 4)$.
A computational check using \textsf{GAP} (or a manual examination of cases) shows that
there is no $C\in J$
such that $(P,Q,A,C)\in \Sigma$ and $\lambda(P,C)=\lambda(Q,C)=()$.
Nonetheless, for $R=\{\{ 1\},\{2,3\},\{4,6\},\{5,7\}\}$,\\
$D=\{1,2,6,7\}$, $E=\{1,2,4,5\}$ we have
$(P,R,A,D)\in\Sigma$ with $\lambda(P,D)=\lambda(R,D)=()$, and $(Q,R,A,E)\in\Sigma$ with
$\lambda(Q,E)=\lambda(R,E)=()$.
\end{runex}

\section{Generators labelled by contiguous cycles}
\label{sec6}

Keeping with the proof outline from Section \ref{sec3}, in this section we deal with
generators $f_{P,A}$ labelled by permutations with descent number $1$.
A typical such permutation is
$$
\xi_{k,l}=(k+l\ \dots\ k+1\ k)\ (1\leq k<k+l\leq r);
$$ 
(see  Lemma \ref{smalldescent}), the notation we will use throughout.
For $l=1$ we get a Coxeter transposition $\xi_{k,k+1}=(k\ k+1)$.

Our aim in this section is two-fold: Firstly, we establish that the generators labelled by 
the same $\xi_{k,l}$ are equal as a consequence of presentation $\present$. This is Lemma \ref{cyclesequal}.
Our second aim is to show that every generator labelled by a cycle $\xi_{k,l}$ of length $l>2$,
can be expessed as a product of a generator labelled by a Coxeter transposition and a generator labelled
by a cycle of length $l-1$; this is Lemma \ref{cycleseliminate1}.
Inductively, this lets us eliminate all the generators labelled by permutations of descent number $1$, except for
one generator labelled by every Coxeter transposition.

\begin{lemma}
\label{smallestA}
Let $k,k+l\in [1,r]$, $l>0$.
The $\lexless$ smallest $A\in J$ for which there exists $P\in I$ such that $A\perp P$ and
$\lambda(P,A)=\xi_{k,l}$ is
$$
A_\mmin =[ 1,k-1]\cup [k+1,r+1].
$$
\end{lemma}

\begin{proof}
Define $P_\mmin =\{P_1,\dots,P_r\}$ as follows:
\begin{eqnarray*}
&& P_i=\{i\} \ (i\in [2,k-1]),\\
&&P_k=\{ k,k+l+1\},\ (P_1=\{ 1,l+2\}\cup [r+2,n] \text{ if } k=1),\\
&&P_i=\{i\} \ (i\in [k+1,k+l]),\\
&&P_i=\{i+1\} \ (i\in [k+l+1,r],\\
&&P_1=P\setminus (P_2\cup\dots\cup P_r)\ (\text{if } k\neq 1).
\end{eqnarray*}

Let us compute $\lambda(P_\mmin ,A_\mmin )$.
(This is our promised sample computation of a label. The computations of labels in subsequent proofs are omitted, but they
can all be done by following exactly the same procedure.)

Clearly,
$$
A(P_\mmin )=[1,k+l]\cup [k+l+2,r+1].
$$
The mapping $\rho_{A(P_\mmin )}\rrestriction_{[1,r]} $ is the unique order preserving bijection
$[1,r]\rightarrow A(P_\mmin )$:
$$
\rho_{A(P_\mmin )}\rrestriction_{[1,r]} \st i\mapsto
\begin{cases}
i, & \text{if } i\in [1,k+l],\\
i+1 & \text{if } i\in [k+l+1,r].
\end{cases}
$$
Likewise, $\rho_{A_\mmin }^{-1}\rrestriction_{A_\mmin }$ is the unique order preserving bijection
$A_\mmin \rightarrow [1,r]$:
\[
\rho_{A_\mmin}^{-1} \rrestriction_{A_\mmin } \st i\mapsto
\begin{cases}
i, & \text{if } i\in [1,k-1],\\
i-1, & i\in [k+1,r+1].
\end{cases}
\]
The final ingredient needed for computing $\lambda(P_\mmin ,A_\mmin )$ is $\gamma_{P_\mmin ,A_\mmin }$,
which  maps $A(P_\mmin )$ into $A_\mmin $ by sending each $\min P_i$  into the unique element of $A_\mmin \cap P_i$:
\[
\gamma_{P_\mmin ,A_\mmin }\st i\mapsto 
\begin{cases}
i, & \text{if } i \in [1,r+1]\setminus\{ k,k+l+1\},\\
k+l+1, & \text{if } i=k.
\end{cases}
\]
So now, for $i\in [1,r]$, we have:
\[
\rho_{A(P_\mmin )} \gamma_{P_\mmin , A_\mmin }\rho_{A_\mmin }^{-1} \st i\mapsto
\begin{cases}
i, & \text{if } i\in [1,k-1],\\
k+l, & \text{if } i=k,\\
i-1, & \text{if } i\in [k+1,k+l],\\
i, & \text{if } i\in [k+l+1,r],
\end{cases}
\]
i.e. $\lambda(P_\mmin ,A_\mmin )=\xi_{k,l}$.

To prove minimality, let $A\in J$ be any set with $A\lexless A_\mmin $, which means that $[1,k]\subseteq A$.
Let $Q\in I$ be any partition with $A\perp Q$.
If $Q=\{ Q_1,\dots, Q_r\}$ with $\min Q_1<\dots <\min Q_r$ it follows that $\min Q_i=i$ for all
$i\in [1,k]$.
But then it follows that $\lambda(Q,A)$ fixes all $i\in [1,k]$, and hence $\lambda(Q,A)\neq \xi_{k,l}$.
\end{proof}

\begin{lemma}
\label{cyclesequal}
Let $k,k+l\in [1,r]$, $l>0$, and 
let $P\in I$, $A\in J$ be such that $A\perp P$ and $\lambda(P,A)=\xi_{k,l}$.
Then the relation $f_{P,A}=f_{P_\mmin ,A_\mmin }$ is a consequence of $\present$,
where $P_\mmin $ and $A_\mmin $ are as in Lemma \ref{smallestA} and its proof.
\end{lemma}

\begin{proof}
Induction on $A$ with respect to $\lexless$.
By Lemma \ref{smallestA}, the induction is anchored for $A=A_\mmin $, in which case the assertion follows by Lemma \ref{lemmav2}.
Suppose $A\lexgreater A_\mmin $, and let $A=\{ a_1,\dots,a_r\}$, $a_1<\dots<a_r$.
Let $t$ be the smallest subscript from $[1,k-1]$ such that $a_t>t$, or, if such does not exist,
the smallest subscript from $[k,r]$ such that $a_t>t+1$.
It follows from $A\lexgreater A_\mmin$ that one of these two cases must arise, and it is clear that
$a_t-1\in [1,r-1]\setminus A$.

Next, let $P=\{ P_1,\dots, P_r\}$, with $p_i=\min P_i$ and $p_1<\dots<p_r$.
Since $\lambda(P,A)=\xi_{k,l}$ we must have
\begin{alignat*}{2}
& a_i\in P_i &\quad& (i\in [1,k-1]),\\
& a_{k+l} \in P_k, &&\\
&a_{i-1}\in P_i && (i\in [k+1,k+l]),\\
&a_i\in P_i && (i\in [k+l+1,r]).
\end{alignat*}

Suppose $a_i\neq p_i$ for some $i\in [1,k-1]$,
and pick the smallest such $i$.
Let $A^\prime=(A\setminus \{a_i\})\cup \{p_i\}$.
Clearly, $A^\prime\perp P$,
and 
$\lambda(P,A^\prime)=\xi_{k,l}$ since by minimality of $i$ we have 
$a_{i-1}=p_{i-1}<p_i<a_i<a_{i+1}$.
By Lemma \ref{lemmav2} the relation $f_{P,A}=f_{P,A^\prime}$ is a consequence of $\present$.
Also, $A^\prime\lexless A$, and so $f_{P,A^\prime}=f_{P_\mmin ,A_\mmin }$ by induction,
implying $f_{P,A}=f_{P_\mmin ,A_\mmin }$ in this case.

So from now on we assume that
\[
a_i=p_i=\min P_i\ (i\in [1,k-1]).
\]
Define a new partition
$Q=\{ Q_1,\dots, Q_r\}$ where
\begin{eqnarray*}
  && Q_i=\{a_i\}\ (i\in [2,k-1]),\\
&& Q_k=\{p_k,a_{k+l}\},\ (Q_1=([1,n]\setminus (Q_2\cup\dots\cup Q_n))\cup\{ a_{l+1}\} \text{ if } k=1),\\
&&Q_i=\{ a_{i-1}\} \ (i\in [k+1,k+l]),\\
&&Q_i=\{a_i\}\ (i\in [k+l+1,r]),\\
&& Q_1=[1,n]\setminus (Q_2\cup\dots\cup Q_r)\ (\text{if } k\neq 1).
\end{eqnarray*}
This partition is obtained from $P$ by reducing each $P_i$ ($i\in [2,r]\setminus\{ k\}$)
to just $P_i\cap A$, reducing $P_k$ to $\{ p_k,a_{k+l}\}$,
and moving all the remaining elements into $Q_1$.
Clearly, we have $A\perp Q$.

Suppose now that $a_t -1\neq p_k$, which means that $a_t -1\in Q_1$.
Define yet another partition $Q^\prime =\{ Q_1^\prime,\dots, Q_r^\prime\}$ by:
\[
Q_1^\prime= Q_1\setminus \{ a_t -1\},\ 
Q_t^\prime = Q_t \cup \{ a_t -1\},\ 
Q_i^\prime=Q_i  \ (i\in [1,r]\setminus \{1,t\}).
\]
Clearly $A\perp Q^\prime$.
Note that $\min Q_i^\prime = \min Q_i$ for all $i\neq t$, and that either $\min Q_t^\prime = \min Q_t$ or else $\min Q_t^\prime= (\min Q_t)-1$.
In any case, we have
$\min Q_1^\prime <\dots<\min Q_r^\prime$, and so $\lambda(Q^\prime,A)=\xi_{k,l}$. 
By  Lemma \ref{lemmav2} this implies $f_{P,A}=f_{Q^\prime,A}$.
Let
\[
A^\prime= (A\setminus \{a_t\})\cup \{ a_t -1\};
\]
again $A^\prime\perp Q^\prime$ and $\lambda(Q^\prime, A^\prime)=\xi_{k,l}$,
so $f_{Q^\prime,A}=f_{Q^\prime,A^\prime}$ by Lemma \ref{lemmav2a}.
But $A^\prime\lexless A$, and so $f_{Q^\prime,A^\prime}=f_{P_\mmin ,A_\mmin }$ by induction,
proving $f_{P,A}=f_{P_\mmin ,A_\mmin }$ in this case.

Finally, consider the case where $a_t -1=p_k$.
Since we have $a_{k-1}=p_{k-1}<p_k<a_k$ it follows that $t=k$, 
and so, by the definition of $t$ from the beginning of the proof, we have $a_k>k+1$.
Since $a_{k-1}=k-1$ it follows that $k\in Q_1$.
Define a new partition $Q^{\prime\prime}=\{ Q_1^{\prime\prime},\dots, Q_r^{\prime\prime}\}$ by:
\begin{eqnarray*}
&Q_1^{\prime\prime} =Q_1\setminus\{k\},\ 
Q_k^{\prime\prime}=\{k,a_{k+l}\},\ 
Q_{k+1}^{\prime\prime}=\{ a_k-1,a_k\},&\\
&Q_i^{\prime\prime}=Q_i \ (i\in [1,r]\setminus\{ 1,k,k+1\}).&
\end{eqnarray*}
Yet again $A\perp Q^{\prime\prime}$ and $\lambda(Q^{\prime\prime},A)=\xi_{k,l}$,
so that $f_{P,A}=f_{Q^{\prime\prime},A}$ by Lemma \ref{lemmav2}.
Next let $A^{\prime\prime}=(A\setminus\{ a_k\})\cup\{ a_k-1\}$.
Once more $A^{\prime\prime}\perp Q^{\prime\prime}$ and
$\lambda(Q^{\prime\prime},A^{\prime\prime})=\xi_{k,l}$.
The relation $f_{Q^{\prime\prime},A^{\prime\prime}}=f_{Q^{\prime\prime},A}$
follows by Lemma \ref{lemmav2a},
while $f_{Q^{\prime\prime},A^{\prime\prime}}=f_{P_\mmin ,A_\mmin }$ follows by induction since $A^{\prime\prime}\lexless A$,
completing the proof of this final case, and hence of the lemma.
\end{proof}

\begin{remark}
\label{disconnected}
The above proof can be interpreted as follows: Suppose that $\pi\in S_r$ is a permutation which appears as a label, i.e.
suppose that the set 
\[
V(\pi)=\{ (P,A)\in I\times J\st \lambda(P,I)=\pi\}
\]
 is non-empty.
Let 
\[
E(\pi)=\{ ((P_1,A_1),(P_2,A_2))\in V(\pi)\times V(\pi) \st P_1=P_2 \text{ or } A_1=A_2\},
\]
and let the graph $G(\pi)$ be $(V(\pi),E(\pi))$.
The proof of Lemma \ref{cyclesequal} ascertains that for $\pi=\xi_{k,l}$ the graph $G(\pi)$ is connected.
It is tempting to hope that this might be the case for an arbitrary $\pi$,
and indeed a quick computational check shows that this is the case in our running example.
Unfortunately, not so in general: for $n=7$, $r=5$ the permutation $\pi=(2,3)(4,5)$ appears exactly twice as the label,
and the corresponding partitions and subsets are $(P_1,A_1)=(\{\{ 1\},\{2,4\},\{3\},\{5,7\},\{6\}\}, \{1,3,4,6,7\})$
and $(P_2,A_2)=(\{\{ 1\},\{2,5\},\{3\},\{4,7\},\{6\}\},\{1,3,5,6,7\})$.
\end{remark}

\begin{lemma}
\label{cycleseliminate}
For any $k,k+l\in [1,r]$, $l\geq 2$, there exists  $(P,Q,A,B) \in \Sigma$ such that:
\begin{alignat}{2}
\label{labelsP}
 & \lambda(P,A) = (), & \quad & 
\lambda(P,B) =\xi_{k+1,l-1},\\
\label{labelsQ}
&\lambda(Q,A) = (k \ k+1), & &
\lambda(Q,B) = \xi_{k,l}.
\end{alignat}
\end{lemma}

\begin{proof}
Define $A$, $B$, $P=\{ P_1,\dots, P_r\}$ and $Q=\{ Q_1,\dots, Q_r\}$ by
\begin{eqnarray*}
&&A=[1,r+2]\setminus \{k,k+l+2\},\\
&&B=[1,r+2]\setminus \{k,k+2\},\\
&&P_1=\{1\}\cup [r+3,n]\ (\text{if } k\neq 1),\\
&& P_i =\{ i \}  \ (i\in [2,  k-1]),\\
&& P_k=\{k, k+1\},  \ (P_1=\{1,2\}\cup [r+3,n] \text{ if } k=1),\\
&& P_{k+1}=\{k+2, k+l+2\}, \\
&&P_i=\{i+1\} \ (i\in [k+2,k+l]), \\
&&P_i=\{i+2\} \ (i\in [k+l+1,r]),\\
&&Q_1=\{1\}\cup [r+3,n]\ (\text{if } k\neq 1),\\
&&Q_i = \{ i \}  \ (i \in [2, k-1]),\\
&&Q_k=\{ k, k+2,k+l+2\}, \ (Q_1=\{1,3,l+3\}\cup [r+3,n] \text{ if } k=1),\\
&&Q_{k+1}=\{ k+1 \},\\
&&Q_i=\{i+1\}  \ (i \in [k+2,k+l]),\\
&&Q_i=\{i+2\}  \ (i\in [k+l+1,r]).
\end{eqnarray*}
It  is easily observed that $\{A,B\}\perp \{ P,Q\}$,
and a routine computation of labels shows that
\eqref{labelsP} and \eqref{labelsQ} hold.
Now note that
\[
\lambda(P,A)^{-1} \lambda(P,B) =
\xi_{k+1,l-1}=\lambda(Q,A)^{-1} \lambda(Q,B),
\]
and hence $(P,Q,A,B)$ is a singular square by Lemma \ref{singularity}.
\end{proof}

\begin{lemma}
\label{cycleseliminate1}
Let $1 \leq k \leq k+l \leq r$ and $l \geq 2$. 
For every $f\in F$ with $\lambda(f)=\xi_{k,l}$ there exist
$f_1,f_2\in F$ such that $\lambda(f_1)=(k\ k+1)$, $\lambda(f_2)=\xi_{k+1,l-1}$
and the relation $f=f_1f_2$ is a consequence of $\present$. 
\end{lemma}

\begin{sloppypar}
\begin{proof}
From the singular square $(P,Q,A,B)$ exhibited in the previous lemma, 
we get the relation
$f_{P,A}^{-1} f_{P,B} = f_{Q,A}^{-1} f_{Q,B}$ in \eqref{eqn_bottom}.
Since $\lambda(P,A)=()$ it follows from Lemma \ref{fixlemma1} that $f_{P,A}=1$. Hence we
have
$f_{Q,B}=f_{Q,A}f_{P,B}$, and the labels are as required.
The assertion now follows for \emph{all} generators with label $\xi_{k,l}$,
as they are all equal to $f_{Q,B}$ by Lemma \ref{cyclesequal}.
\end{proof}
\end{sloppypar}

\begin{remark}
\label{notinplace}
It is not true that every generator $f_{Q,B}$ with label $\xi_{k,l}$ ($l\geq 2$) can be eliminated \emph{in place}, i.e. by means of a
singular square $(P,Q,A,B)$ as in the above lemma.
For instance, in our Running Example, take $Q=\{ \{1,4,5,7\},\{2\},\{3\},\{6\}\}$ and $B=\{2,3,6,7\}$.
Computational search using \textsf{GAP} shows that there do not exist $Q\in I$ and $B\in J$ such that
$(P,Q,A,B)\in\Sigma$ and the four labels are $()$, $(3\ 2\ 1)$, $(1\ 2)$, $(4\ 3\ 2\ 1)$.
In fact there are unique $P$ and $A$ (namely $P=\{\{1,3\},\{2\},\{4,6\},\{5,7\}\}$ and $A=\{2,3,5,6\}$)
such that $(P,Q,A,B)\in\Sigma$, $\lambda(P,B)=(3\ 2\ 1)$ and $\lambda(Q,A)=(1\ 2)$.
But $\lambda(Q,B)=(1\ 2)(3\ 4)$, a permutation of descent number $2$!
\end{remark}

\section{Generators labelled by permutations with descent number \boldmath{$>1$}}
\label{sec7}
The results from Sections \ref{sec4} and \ref{sec6} mean that
 all generators whose labels  have descent number $\leq 1$ can be expressed as products of generators labelled by  Coxeter transpositions, and thus eliminated from $\present$.
Furthermore, we know that every two generators labelled by the same Coxeter transposition are equal.
In this section we show how to eliminate generators whose labels have descent number greater than $1$.
Before we state and prove the main result of the section, let us recall the `image' notation for permutations:
$[l_1,\dots, l_r]$ stands for the permutation $\pi\in S_r$ such that
$i\pi = l_i$ ($i\in [1,r]$).

\begin{lemma}
\label{lem_elimination}
Let $P\in I$, $A\in J$ be such that $A\perp P$ and 
\[
\lambda(P,A)=[l_1,\dots, l_r],
\]
with $\Delta(\lambda(P,A))=d>1$.
Let $v,v+w\in [1,r]$ be such that $l_v$ is the rightmost descent start in $\lambda(P,A)$,
and $(l_v,l_{v+w})$ is the rightmost descent starting at $l_v$.
Then there exist $Q\in I$, $B\in J$ such that
\begin{enumerate}[label=\textup{(\roman*)}, leftmargin=*,widest=2]
\item
\label{lem_eliminationi}
$(P,Q,A,B)\in \Sigma$;
\item
\label{lem_eliminationii}
$\lambda(P,B)=(v+w,\dots,v+1,v)$, a permutation with descent number $1$;
\item
\label{lem_eliminationiii}
\begin{sloppypar}
$\lambda(Q,A)=[l_1,\dots,l_{v-1},l_{v+1},\dots,l_{v+w},l_v,l_{v+w+1},\dots,l_r]$,
a permutation with descent number $d-1$;
\end{sloppypar}
\item
\label{lem_eliminationiv}
The relation $f_{P,A}=f_{P,B}f_{Q,A}$ is a consequence of $\present$.
\end{enumerate}
\end{lemma}

\begin{proof}
Suppose
\begin{eqnarray*}
&&A=\{ a_1,\dots,a_r\},\ a_1<\dots<a_r,\\
&& P=\{P_1,\dots,P_r\},\ p_i=\min P_i,\ p_1<\dots<p_r.
\end{eqnarray*}
The assumption $\lambda(P,A)=[l_1,\dots,l_r]$ means that we have $a_{l_i}\in P_i$.
Let
\[
B=\{ p_1,\dots, p_{v-1},a_{l_{v+1}},\dots, a_{l_{v+w}}, a_{l_v},a_{l_{v+w+1}},\dots,a_{l_r}\};
\] 
note that, since $l_{v+w}<l_v <l_{v+w+1}$, we have
\[
p_1<\dots <p_{v-1} <a_{l_{v+1}}<\dots < a_{l_{v+w}} < a_{l_v}< a_{l_{v+w+1}}<\dots <a_{l_r}.
\]
Also define $Q=\{Q_1,\dots,Q_r\}$ by:
\begin{eqnarray*}
&&Q_i=P_i\ (i\in [2,v-1]),\\
&&Q_i=\{ a_{l_{i+1}}\} \ (i\in [v,v+w-1]),\ (Q_1=([1,n]\setminus A)\cup \{ a_{l_2} \} \text{ if } v=1),\\
&& Q_{v+w}=\{ a_{l_v}\},\\
&& Q_i=\{ a_{l_i}\}\ (i\in [v+w+1,r]),\\
&& Q_1=[1,n]\setminus (Q_2\cup\dots\cup Q_r)\ (\text{if } v\neq 1).
\end{eqnarray*}
Essentially, $Q$ is obtained from $P$ by retaining $P_2,\dots,P_{v-1}$,
reducing each $P_v,\dots, P_r$ to its one element intersection with $A$, and moving all the other elements to $P_1$.

It is easy to see that $\{A,B\}\perp \{P,Q\}$, and a routine computation 
shows that \ref{lem_eliminationii}, \ref{lem_eliminationiii} hold, that $\lambda(Q,B)=()$, and
and that $(P,Q,A,B)$ is a singular square by virtue of satisfying condition \ref{SQ3} of Lemma \ref{singularity}.
Thus the relation $f_{P,A}^{-1}f_{P,B}=f_{Q,A}^{-1}f_{Q,B}$ is in \eqref{eqn_bottom}.
By Lemma \ref{fixlemma1} we have $f_{Q,B}=1$, and so we obtain $f_{P,A}=f_{P,B}f_{Q,A}$,
as required.
\end{proof}

\begin{runex}
For $P=\{ \{1,7\},\{2,5\},\{3,6\},\{4\}\}$ and $A=\{4,5,6,7\}$ we have
$\lambda(P,A)=[4,2,3,1]$, a permutation of descent number $3$.
Following the above proof, we define
$B=\{ 1,2,4,6\}$ and $Q=\{ \{ 1,3,7\}, \{2,5\},$ $\{4\},\{6\}\}$.
We have $(P,Q,A,B)\in \Sigma$, and the labels are
$\lambda(P,B)=(4\ 3)$  (descent number $1$),
$\lambda(Q,A)=[4,2,1,3]$ (descent number $2$) and
$\lambda(Q,B)=()$.
\end{runex}

\begin{remark}
In contrast with the observation made in Remark \ref{notinplace}, note that 
Lemma \ref{lem_elimination} allows us to eliminate every generator $f_{P,A}$ with
$\Delta(\lambda(P,A))>1$ \emph{in place where it occurs}.
This is important as we have already indicated (Remark \ref{disconnected}) that we cannot prove equality of all the generators 
with equal labels of descent number $>1$ just by linking them via singular squares having two adjacent vertices labelled $()$.
This fact will also be a key technicality working in the background of the proof of Lemma \ref{coxrel3} below.
\end{remark}

\section{Coxeter relations}
\label{sec8}

At this stage of the proof we know that all the generators from presentation $\present$ that are not labelled by Coxeter transpositions are redundant, and that any two generators labelled by the same Coxeter transposition are equal.
There remains to be proved that the latter generators satisfy 
the standard Coxeter relations (see \eqref{Cox} in Section \ref{sec3}).

\begin{lemma}
\label{coxrel1}
For any $k\in [1,r-1]$ and any $f\in F$ with $\lambda(f)=(k\ k+1)$ the relation $f^2=1$ is a consequence of $\present$. 
\end{lemma}

\begin{proof}
By Lemma \ref{cyclesequal}, it suffices to find a single such $f$.
Define $A,B\in J$, $P=\{ P_1,\dots, P_r\}\in I$, $Q=\{ Q_1,\dots, Q_r\}\in I$ by:
\begin{eqnarray*}
&&A = [1,r+2] \setminus \{ k, k+3 \},
 \\
&& B = [1,r+2] \setminus \{ k, k+1 \}, \\
&& P_1=\{1 \} \cup [r+3,n] \ (\text{if } k\neq 1),\\ 
&& P_i=\{ i \}   (i \in [2,k-1]),  \\
&&P_k=\{ k,k+2 \},\ (P_1=\{1,3\}\cup [r+3,n] \text{ if } k=1),\\
&&P_{k+1}=\{ k+1, k+3 \}, \\
&&P_i=\{i+2\}   \ (i \in [k+2,r]),\\ 
&& Q_1=\{1, k\} \cup [r+3,n] \ (\text{if } k\neq 1),\\
&&Q_i=\{ i \}   \ (i \in [2,k-1]),  \\
&&Q_k=\{k+1,k+3\},  \ (Q_k=\{1,2,4\} \cup [r+3,n] \text{ if } k=1), \\
&&Q_{k+1}=\{k+2\},   \\
&&Q_i=\{i+2\}   \ (i \in [k+2,r]).
\end{eqnarray*}
It is clear that $\{A,B\}\perp\{P,Q\}$.
A routine label computation shows that
\[
\lambda(P,A)=\lambda(Q,B)=(k\ k+1),\ 
\lambda(P,B)=\lambda(Q,A)=().
\]
By Lemma \ref{singularity}, the square $(P,Q,A,B)$ is singular.
It yields the relation $f_{P,A}^{-1}f_{P,B}=f_{Q,A}^{-1}f_{Q,B}$.
By Lemmas \ref{fixlemma1} and  \ref{cyclesequal} we have $f_{P,B}=f_{Q,A}=1$ and $f_{P,A}=f_{Q,B}$.
Therefore $f_{P,A}^2=1$, as required.
\end{proof}

\begin{runex}
The square demonstrating that $f^2=1$ for any generator with label $(2\ 3)$ in our Running Example is
given by $A=\{1,3,4,6\}$, $B=\{1,4,5,6\}$, $P=\{\{1,7\},\{2,4\},\{3,5\},\{6\}\}$,
$Q=\{\{1,2,7\},\{3,5\},\{4\},\{6\}\}$.
The labels are, as expected, $\lambda(P,A)=\lambda(Q,B)=(2\ 3)$ and
$\lambda(P,B)=\lambda(Q,A)=()$.
\end{runex}

\begin{lemma}
\label{coxrel2}
Let $k,l\in [1,r]$ be such that $k+1<l$.
For any $f,g\in F$ with $\lambda(f)=(k\ k+1)$, $\lambda(g)=(l\ l+1)$ the relation $fg=gf$ is a consequence of $\present$.
\end{lemma}

\begin{proof}
By Lemma \ref{cyclesequal} it suffices to prove the assertion for a particular pair $f,g$.
Let us define $A,B,C\in J$, $P=\{ P_1,\dots,P_r\}\in I$, $Q=\{ Q_1,\dots, Q_r\}\in I$,
$R=\{ R_1,\dots, R_r\}\in I$ as follows:
\begin{eqnarray*}
&&A  =  [1,r+2] \setminus \{ k+2, l+1 \}, \\
&&B = [1,r+2] \setminus \{ k, l+1 \}, \\
&&C  =  [1,r+2] \setminus \{ k, l+3 \}, \\
&&P_1=\{ 1,l+1 \}  \cup [r+3,n]\ (\text{if } k\neq 1), \\
&& P_i=\{ i\} \ (i \in [2,k-1]), \\
&& P_k=\{k,k+2\},\ (P_1=\{ 1,3,l+1\}\cup [r+3,n] \text{ if } k=1)\\
&&P_{k+1}=\{ k+1 \} , \\
&&P_i=\{i+1\}  \ (i \in [k+2,l-1]), \\
&&P_i=\{i+2\}  \ (i \in [l,r]),\\
&&Q_1=\{ 1 \} \cup [r+3,n] (\text{if } k\neq 1),\\
&&Q_i=\{ i\}  \ (i \in [2,k-1]),\\
&&Q_k=\{k,k+2\} ,\  (Q_1=\{ 1,3\}\cup [r+3,n] \text{ if } k=1)\\
&&Q_{k+1}=\{ k+1 \},\\
&&Q_i=\{i+1\}\ (i \in [k+2,l-1]), \\
&&Q_l=\{l+1,l+3\} ,\\
&&Q_{l+1}=\{l+2\} , \\
&&Q_i=\{i+2\} \ (i \in [l+2,r]),\\
&&R_1=\{ 1 ,k \} \cup [r+3,n] \ (\text{if } k\neq 1), \\
&&R_i=\{ i\}  \ (i \in [2,k-1]),\\
&&R_i=\{i+1\}  \ (i \in [k,l-1]), (R_1=\{ 1,2\}\cup [r+3,n] \text{ if } k=1),\\
&&R_l=\{ l+1,l+3 \}, \\
&&R_{l+1}=\{l+2\}  , \\
&&R_i=\{i+2\}  \ (i \in [l+2,r]).
\end{eqnarray*}
It is easy to see that $A\perp \{P,Q\}$, $B\perp\{P,Q,R\}$, $C\perp\{Q,R\}$.
The labels are:
\begin{align*}
\lambda(P,A) &= (), & \lambda(P,B) &= (k\ k+1), & & \\
\lambda(Q,A) &= (l\ l+1), & \lambda(Q,B) &=(k\ k+1)(l\ l+1), & \lambda(Q,C) &=(k\ k+1),\\
& & \lambda(R,B) &=(l,l+1), & \lambda(R,C) &= ().
\end{align*}
By Lemma \ref{singularity}, we have $(P,Q,A,B),(Q,R,B,C)\in \Sigma$, yielding
$f_{P,A}^{-1}f_{P,B}=f_{Q,A}^{-1}f_{Q,B}$ and $f_{Q,B}^{-1}f_{Q,C}=f_{R,B}^{-1} f_{R,C}$.
Eliminating $f_{Q,B}$ gives
$f_{P,B}^{-1}f_{P,A}f_{Q,A}^{-1}=f_{R,B}^{-1}f_{R,C}f_{Q,C}^{-1}$.
By Lemmas \ref{fixlemma1} and  \ref{cyclesequal} we can eliminate $f_{P,A}=f_{R,C}=1$, $f_{Q,C}=f_{P,B}$ and $f_{R,B}=f_{Q,A}$,
giving us $f_{P,B}f_{Q,A}=f_{Q,A}f_{P,B}$, as required.
\end{proof}

\begin{runex}
In order to exhibit the relation $fg=gf$ with $\lambda(f)=(1\ 2)$, $\lambda(g)=(3\ 4)$ in our Running Example, one ought to take:
$A=\{1,2,5,6\}$, $B=\{2,3,5,6\}$, $C=\{2,3,4,5\}$, $P=\{\{1,3,4,7\},\{2\},\{5\},\{6\}\}$,
$Q=\{\{1,3,7\},\{2\},\{4,6\},\{5\}\}$,
$R=\{\{1,2,7\},\{3\},\{4,6\},\{5\}\}$.
\end{runex}

\begin{lemma}
\label{coxrel3}
For any $k\in [1,r-2]$ and any $g,h\in F$ with $\lambda(f)=(k\ k+1)$, $\lambda(g)=(k+1\ k+2)$,
the relation $ghg=hgh$ is a consequence of $\present$.
\end{lemma}

\begin{proof}
Again, by Lemma \ref{cyclesequal}, it suffices to prove the assertion for a particular pair $f,g$.
Define $A,B\in J$, $P=\{P_1,\dots,P_r\}\in I$, $Q=\{Q_1,\dots,Q_r\}\in I$ by:
\begin{eqnarray*}
&&A=[1,r+2]\setminus\{k+1,k+4\}
,\\
&&B=[1,r+2]\setminus\{k,k+1\},\\
&&P_1=\{1\}\cup [r+3,n] \ (\text{if } k\neq 1),\\
&&P_i=\{i\}\ (i\in [2,k-1]),\\
&&P_k=\{k,k+1,k+4\}\ (P_1=\{1,2,5\} \cup [r+3,n] \text{ if } k=1),\\
&&P_{k+1}=\{k+2\},\\
&&P_{k+2}=\{k+3\},\\
&&P_i=\{i+2\}\ (i\in [k+3,r],\\
&&Q_1=\{ 1\} \cup [r+3,n] \ (\text{if } k\neq 1),\\
&&Q_i=\{i\}\ (i\in [2,k-1]),\\
&&Q_k=\{k,k+4\},\ (Q_1=\{1,5\} \cup [r+3,n] \text{ if } k=1),\\
&&Q_{k+1}=\{ k+1,k+3\},\\
&&Q_{k+2}=\{k+2\},\\
&&Q_i=\{i+2\}\ (i\in [k+3,r]).
\end{eqnarray*}
Routinely, $\{A,B\}\perp \{P,Q\}$ and
\begin{alignat*}{2}
\lambda(P,A) &= (), &\quad \lambda(P,B) &= (k+2\ k+1\ k),\\
\lambda(Q,A) &= (k+1\ k+2), & \lambda(Q,B)&=(k\ k+2).
\end{alignat*}
Hence $(P,Q,A,B)\in\Sigma$ (Lemma \ref{singularity}), and we have the relation
$f_{P,A}^{-1} f_{P,B}=f_{Q,A}^{-1} f_{Q,B}$
in \eqref{eqn_bottom}.
By Lemma \ref{fixlemma1} we have $f_{P,A}=1$, and so
\begin{equation}
\label{x2}
 f_{P,B}=f_{Q,A}^{-1} f_{Q,B}.
\end{equation}
Note that $\Delta(f_{Q,B})=\Delta((k\ k+2))=2$; Lemma \ref{lem_elimination} applied to $f_{Q,B}$ yields 
$f_{Q,B}=f_1f_2$
for some $f_1,f_2\in F$ with $\lambda(f_1)=(k+2\ k+1\ k)$, $\lambda(f_2)=(k\ k+1)$.
By Lemma \ref{cyclesequal} we have $f_1=f_{P,B}$, and so \eqref{x2} becomes
\begin{equation}
\label{x3}
f_{Q,A}f_{P,B}=f_{P,B} f_2.
\end{equation}
By Lemma \ref{cycleseliminate1} applied to $f_{P,B}$ we have $f_{P,B}=f_3f_4$ for some $f_3,f_4\in F$ satisfying
$\lambda(f_3)=(k\ k+1)$, $\lambda(f_4)=(k+1\ k+2)$.
By Lemma \ref{cyclesequal} we have $f_3=f_2$ and $f_4=f_{Q,A}$.
Letting $g=f_2$ and $h=f_{Q,A}$, and substituting into \eqref{x3} we obtain $ghg=hgh$, as required.
\end{proof}

\begin{runex}
For $k=2$ the singular square constructed above has
$A=\{ 1,2,4,5\}$, $B = \{1,4,5,6\}$, 
$P=\{\{ 1,7\},\{2,3,6\},\{4\},\{5\}\}$ and
$Q=\{\{1,7\},\{2,6\},\{3,5\},\{4\}\}$.
\end{runex}

With the proof of Lemma \ref{coxrel3}, the prof of our Main Theorem, as outlined in Section \ref{sec3}, is complete.

\section{Concluding remarks}

Obviously, one would quite like to be able to describe completely the structure of the free idempotent generated semigroup
$\IG(E(T_n))$, the Main Theorem providing an essential ingredient of such a description.
In fact, the Main Theorem does enable such a description of a close relative:
For a general \emph{regular} semigroup $S$, with the set of idempotents $E=E(S)$, the 
\emph{free regular idempotent generated semigroup} $\RIG(E)$ is the free object in the 
the category of regular idempotent generated semigroups with biordered set of idempotents isomorphic to $E$.
For any two idempotents $e,f\in E$ their \emph{sandwich set} is defined by
$S(e,f)=\{ h\in E\st ehf=ef,\ fhe=h\}$, and is known to be non-empty.
The semigroup $\RIG(E)$ can be defined as a quotient of $\IG(E)$ by adding the relations
$ehf=ef$ for all $e,f\in E$ and all $h\in S(e,f)$.
For more details see \cite{nambooripad80} or \cite{margolismeakin09}.
In particular, as pointed out in \cite[Theorem 3.6]{margolismeakin09}, the maximal subgroups of $\IG(E)$ and $\RIG(E)$
coincide. So our Main Theorem remains valid verbatim if $\IG(E)$ is replaced by $\RIG(E)$.

It is well known that $T_n$ decomposes into a chain of $\GD$-classes $D_n$, $D_{n-1}$, \dots , $D_1$,
where $D_r$ consists of all mappings of rank $r$. 
Then $\RIG(E(T_n))$ also has a chain of $\GD$-classes $\overline{D}_n$, $\overline{D}_{n-1}$,
$D_{n-2}$, \ldots, $D_1$.
Here $\overline{D}_n$ consists of a single element, an indecomposable identity for the whole semigroup.
The next $\GD$-class $\overline{D}_{n-1}$ is obtained from $D_{n-1}$ by replacing the underlying group $S_{n-1}$
with the free group $F$ of rank $\binom{n}{2}-1$.
The structure matrix for the principal ideal is obtained by taking the matrix for $D_{n-1}$ and replacing 
each non-zero entry indexed by $(P,A)$ by the generator $f_{P,A}$ of $F$.
The remaining $\GD$-classes of $\RIG(E(T_n))$ are exactly the same as those of $T_n$.
The products within $\overline{D}_{n-1}$ are governed by its Rees matrix representation.
All the remaining products are exactly the same as in $T_n$, of course with the elements of $\overline{D}_{n-1}$ replaced by the corresponding elements of $D_{n-1}$ via the obvious natural homomorphism.
All the above claims follow easily from the Main Theorem, presentation $\present$ exhibited in Section \ref{sec3}
applied to $D_{n-1}$ upon noting that there are no singular squares in this case, and the basic properties
of free idempotent generated semigroups as listed in \cite[Section 1]{grayta}.

The above remarks underline the timeliness of shifting the focus of investigations in this area from maximal subgroups to enhancing our understanding of the general structure of $\IG(E)$.

The \emph{dual} of a semigroup $S$ is the semigroup $S^{\text{op}}$ with the same underlying set and multiplication
$\ast$ defined by $x\ast y = yx$. 
In the dual of the full transformation semigroup $T_n$ the mappings are written to the left of their arguments and are composed from right to left.
From the definition of the free idempotent generated semigroup \eqref{eq1} it is obvious that
$\IG(E(S^{\text{op}}))\cong (\IG(E(S))^{\text{op}}$.
Since every group is isomorphic to its dual (because of the anti-isomorphism $x\mapsto x^{-1}$) it follows
that maximal subgroups of $S$ and $S^{\text{op}}$ coincide.
In particular, our Main Theorem remains valid if $T_n$ is replaced by its dual.

\begin{flushleft}
Centro de \'{A}lgebra da Universidade de Lisboa \\ 
Av. Prof. Gama Pinto 2  \\
1649-003 Lisboa,  Portugal.\\
\smallskip
\texttt{rdgray@fc.ul.pt}
\end{flushleft}

\begin{flushleft}
School of Mathematics and Statistics\\
University of St Andrews\\
St Andrews KY16 9SS\\
Scotland, U.K.\\
\smallskip
\texttt{nik@mcs.st-and.ac.uk}
\end{flushleft}

\end{document}